\documentclass[11pt]{article}
\usepackage{cite}
\usepackage{mathrsfs}
\usepackage{amssymb,amsfonts,amsmath}
\usepackage{enumerate}
\usepackage{indentfirst}
\usepackage{latexsym,bm}
\usepackage[noend]{algpseudocode}
\usepackage{algorithmicx,algorithm}
\usepackage{algorithm}
\usepackage{makeidx}
\usepackage{fancybox}
\usepackage{color}
\usepackage{multicol}
\usepackage[latin1]{inputenc}
\usepackage{graphicx}
\usepackage{epstopdf}
\usepackage{pstricks,pst-node,pst-text,pst-3d}
\usepackage{tikz}
\newtheorem{thm}{Theorem}[section]

\newtheorem{lem}[thm]{Lemma}

\newtheorem{obs}[thm]{Observation}

\newenvironment{pf}{{\noindent \it \bf Proof:}}{{\hfill$\Box$}\\}

\def\qed{\hfill \nopagebreak\rule{5pt}{8pt}}

\begin{document}

\title{\bf Semicomplete Compositions of Digraphs}
\author{Yuefang Sun\footnote{yuefangsun2013@163.com.}{ }\\
%\footnote{The author was supported by Zhejiang Provincial Natural Science Foundation (No. LY20A010013) and National Natural Science Foundation of China (No. 11401389).}{ }
%School of Mathematics and Statistics, Ningbo University,\\ Zhejiang 315211, P. R. China\\
Department of Mathematics, Shaoxing University\\ Zhejiang 312000, P. R. China}
%Department of Mathematics, Shaoxing University\\ Zhejiang 312000, P. R. China, }
%$^{2}$ Department of Computer Science\\
%Royal Holloway, University of London\\
%Egham, Surrey, TW20 0EX, UK\\
%g.gutin@rhul.ac.uk, ai18810910153@163.com}
%\date{}
\maketitle

\begin{abstract}
Let $T$ be a digraph with vertices $u_1, \dots, u_t$ ($t\ge 2$) and let $H_1, \dots, H_t$ be digraphs such that $H_i$ has vertices $u_{i,j_i},\ 1\le j_i\le n_i.$
Then the composition $Q=T[H_1, \dots, H_t]$ is a digraph with vertex set $\{u_{i,j_i}\colon\, 1\le i\le t, 1\le j_i\le n_i\}$ and arc set
$$A(Q)=\cup^t_{i=1}A(H_i)\cup \{u_{ij_i}u_{pq_p}\colon\, u_iu_p\in A(T), 1\le j_i\le n_i, 1\le q_p\le n_p\}.$$ The composition $Q=T[H_1, \dots, H_t]$ is a semicomplete composition if $T$ is semicomplete, i.e. there is at least one arc between every pair of vertices. Digraph compositions generalize some families of digraphs, including (extended) semicomplete digraphs, quasi-transitive digraphs and lexicographic product digraphs. In particular, strong semicomplete compositions form a significant generalization of strong quasi-transitive digraphs.

In this paper, we study the structural properties of semicomplete compositions and obtain results on connectivity, paths, cycles, strong spanning subdigraphs and acyclic spanning subgraphs. Our results show that this class of digraphs shares some nice properties of quasi-transitive digraphs. %Several open questions are also posed.
\vspace{0.3cm}
\vspace{0.3cm}

\end{abstract}

\section{Introduction}

We refer the readers to \cite{Bang-Jensen-Gutin, Bondy} for graph-theoretical notation and terminology not given here.  We use $[n]$ to denote the set of all natural numbers from 1 to $n$. Let $\overline{K}_p$ stand for the digraph of order $p$ with no arcs, and $\overrightarrow{C}_k$ and $\overrightarrow{P}_k$ denote the cycle and path with $k$ vertices, respectively. We call a vertex $x$ in a digraph a {\em sink} (resp. {\em source}) if $d^+(x)=0$ (resp. $d^-(x)=0$).

Let $D$ be a digraph. If there is an arc from a vertex $x$ to a vertex $y$ in $D$, then we say that $x$ {\em dominates} $y$ and denote it by $x\rightarrow y$. If there is no arc from $x$ to $y$ we shall use the notation $x\not\rightarrow y$. If $A$ and $B$ are two subdigraphs of $D$ and every vertex of $A$ dominates each vertex of $B$, then we say that $A$ {\em dominates} $B$ and denote it by $A\rightarrow B$. We shall use $A\Rightarrow B$ to denote that $A$ dominates $B$ and there is no arcs from $B$ to $A$. A digraph $D$ is {\em triangular} with a partition $\{V_0, V_1, V_2\}$, if $V(D)$ can be partitioned into three disjoint sets $V_0, V_1, V_2$ with $V_0\Rightarrow V_1\Rightarrow V_2\Rightarrow V_0$.

Let $\{H_u\mid u\in V(D)\}$ be a disjoint collection of digraphs indexed by the vertices of $D$. To {\em substitute} $H_u$ for every $u$ means to obtain a new digraph $D'$ from $D$ by replacing $u$ with $H_u$ such that $H_u\rightarrow H_v$ in $D'$ if and only if $u\rightarrow v$ in $D$.
Let $H$ be a subdigraph of $D$. To {\em shrink} $H$ in $D$ means to obtain a new digraph $D'$ from $D$ by replacing $H$ with a vertex $v$ so that in $D'$ a vertex $x \not\in V(H)$ dominates (resp. is dominated by) $v$ if and only if $D$ contains an arc from $x$ to $H$ (resp. from $H$ to $x$).

A digraph $D$ is {\em connected} if $U(D)$ is connected, where $U(D)$ is the {\em underlying graph} of $D$. The {\em complementary graph} of a graph $G$ is denoted by $\overline{G}$. A {\em strong component} of a digraph $D$ is a maximal subset $S\subseteq V(D)$ such that $d_D(x,y)$ is finite for every pair of distinct vertices $x, y\in S$, where $d_D(x, y)$ denotes the length of a shortest path from $x$ to $y$ in $D$ (if there is no such path $d_D(x, y)$ is infinite). An {\em initial} (resp. {\em terminal}) {\em strong component} is one which has no arcs entering (resp. leaving) it. We call a set $S\subseteq V(D)$ a {\em separator} of $D$ if $D-S$ is not strong. A separator $S$ is {\em minimal} if no proper subset of $S$ is a separator. A digraph is {\em $k$-strong-connected} if every separator has at least $k$ vertices. A {\em minimum separator} is a minimum cardinality separator. Clearly, a minimum separator is also a minimal separator. A set $S$ of vertices in a digraph $D$ is {\em independent} if $D[S]$ has no arcs.

A digraph is {\em acyclic} if it has no dicycle. An {\em out-tree} (resp. {\em in-tree}) {\em rooted at a vertex $r$} is an orientation of a tree such that the in-degree (resp. out-degree) of every vertex but $r$ equals one. An {\em out-branching} $B^+_r$ (resp. {\em in-branching} $B^-_r$) in a digraph $D$ is a spanning subdigraph of $D$ which is out-tree (resp. in-tree).

A digraph $D$ is {\em pancyclic} if $D$ contains a cycle of length $k$ for each $3\leq k\leq n$, and is {\em vertex-pancyclic} if every vertex of $D$ is contained in a cycle of length $k$ for each $3\leq k\leq n$. A cycle (path) of a digraph $D$ is {\em Hamiltonian} if it contains all the vertices of $D$. A digraph is {\em Hamiltonian} if it has a Hamiltonian cycle. A {\em $k$-path-cycle subdigraph} $\mathcal{F}$ of a digraph $D$ is a collection of $k$ paths $P_1, \dots, P_k$ and $t$ cycles $C_1, \dots, C_t$ such that all of $P_1, \dots, P_k, C_1, \dots, C_t$ are pairwise disjoint (possibly, $k=0$ or $t=0$). We will denote $\mathcal{F}$ by $\mathcal{F}=P_1\cup \dots \cup P_k\cup C_1\cup \dots \cup C_t$. Furthermore, a {\em $k$-path-cycle factor} is a spanning $k$-path-cycle subdigraph. If $t=0$, $\mathcal{F}$ is a {\em $k$-path subdigraph} and it is a {\em $k$-path factor} (or just a {\em path-factor}) if it is spanning. If $k=0$, we say that $\mathcal{F}$ is a {\em $t$-cycle subdigraph} (or just a {\em cycle subdigraph}) and it is a {\em $t$-cycle factor} (or just a {\em cycle factor}) if it is spanning. The {\em path covering number} of a digraph $D$, denoted $pc(D)$, is the smallest $k$ for which $D$ has a {\em $k$-path factor}. Note that a Hamiltonian path is a 1-path factor.

%A collection $F$ of pairwise vertex disjoint paths and cycles of a digraph $D$ is called a {\em $k$-path-cycle factor} of $D$ if $F$ covers $V(D)$ and has exactly $k\geq 0$ paths. $F$ is called a {\em $k$-path factor} if it contains only paths. We shall call a 0-path-cycle factor a {\em cycle factor}. A {\em cycle subgraph} is a collection of vertex disjoint cycles.

%A {\em $k$-king} in a digraph $D$ is a vertex which can reach every other vertex by a directed path of length at most $k$. A 2-king is called a {\em king}, and a {\em non-king} is a vertex which is not a 3-king.

%An {\em out-tree} ({\em in-tree}, respectively) {\em rooted at a vertex $r$} is an orientation of a tree such that the in-degree (out-degree, respectively) of every vertex but $r$ equals one. An {\em out-branching} ( {\em in-branching}, respectively) in a digraph $D$ is a spanning subgraph of $D$ which is out-tree (in-tree, respectively).

A digraph $D$ is {\em semicomplete} if for every pair $x,y$ of distinct vertices of $D$, there is at least one arc between $x$ and $y.$ In particular, a {\em tournament} is a semicomplete digraph without 2-cycle. %A digraph $D$ is {\em locally semicomplete} if the out-neighborhood and in-neighborhood of every vertex of $D$ induce semicomplete digraphs. % A digraph is {\em semicomplete multipartite} if it is obtained from a complete multipartite graph by replacing every edge by an arc or a pair of opposite arcs. A {\em semicomplete tournament} is a semicomplete multipartite digraph with no 2-cycle.
Let $T$ be a digraph with vertices $u_1, \dots, u_t$ ($t\ge 2$) and let $H_1, \dots, H_t$ be digraphs such that $H_i$ has vertices $u_{i,j_i},\ 1\le j_i\le n_i.$ Then the {\em composition} $Q=T[H_1, \dots, H_t]$ is a digraph with vertex set $\{u_{i,j_i}\colon\, 1\le i\le t, 1\le j_i\le n_i\}$ and arc set $$A(Q)=\cup^t_{i=1}A(H_i)\cup \{u_{ij_i}u_{pq_p}\colon\, u_iu_p\in A(T), 1\le j_i\le n_i, 1\le q_p\le n_p\}.$$ The composition $Q=T[H_1, \dots, H_t]$ is a {\em semicomplete composition} if $T$ is semicomplete. If $Q=T[H_1, \dots, H_t]$ and none of the digraphs $H_1, \dots, H_t$ has an arc, then $Q$ is an {\em extension} of $T$. For any set $\Phi$ of digraphs, $\Phi^{ext}$ denotes the (infinite) set of all extensions of digraphs in $\Phi$, which are called {\em extended $\Phi$-digraphs}. Two vertices $x$ and $y$ in an extended semicomplete digraph are said to be {\em similar} if $x, y\in V(H_i)$ for some $i\in [t]$.%By definition, the class of (extended) semicomplete digraphs is a subclass of semicomplete compositions.%If $Q=T[H_1,\dots , H_t]$ and none of the digraphs $H_1,\dots, H_t$ has an arc, then $Q$ is an {\em extension} of $T$.

The following theorem by Bang-Jensen and Huang gives a complete characterization of quasi-transitive digraphs, and the decomposition below is called the {\em canonical decomposition} of a quasi-transitive digraph.

\begin{thm}\label{intro01}\cite{Bang-Jensen-Huang1}
Let $D$ be a quasi-transitive digraph. Then the following assertions hold:
\begin{description}
\item[(a)] If $D$ is not strong, then there exists a transitive oriented
graph $T$ with vertices $\{u_i\mid i\in [t]\}$ and strong quasi-transitive digraphs $H_1, H_2, \dots, H_t$ such that $D = T[H_1, H_2, \dots, H_t]$, where $H_i$ is substituted for $u_i$, $i\in [t]$.
\item[(b)] If $D$ is strong, then there exists a strong semicomplete
digraph $S$ with vertices $\{v_j\mid j\in [s]\}$ and quasi-transitive digraphs $Q_1, Q_2, \dots, Q_s$ such that $Q_j$ is either a vertex or is non-strong and $D = S[Q_1, Q_2, \dots, Q_s]$, where $Q_j$ is substituted for $v_j$, $j\in [s]$.
\end{description}
\end{thm}

Digraph compositions generalize some families of digraphs, including (extended) semicomplete digraphs, quasi-transitive digraphs (by Theorem~\ref{intro01}) and lexicographic product digraphs (when $H_i$ is the same digraph $H$ for every $i\in [t]$, $Q$ is the lexicographic product of $T$ and $H$, see, e.g., \cite{Hammack}). In particular, strong semicomplete compositions generalize strong quasi-transitive digraphs. To see that strong semicomplete compositions form a significant generalization of strong quasi-transitive digraphs, observe that the Hamiltonian cycle problem is polynomial-time solvable for quasi-transitive digraphs \cite{Gutin4}, but NP-complete for strong semicomplete compositions (see, e.g., \cite{Bang-Jensen-Gutin-Yeo}). While digraph composition has been used since 1990s to study quasi-transitive digraphs and their generalizations, see, e.g., \cite{Bang-Jensen-Huang1, GSHC}, the study of digraph decompositions in their own right was initiated only recently by Sun, Gutin and Ai in \cite{Sun-Gutin-Ai}, where some results on the existence of arc-disjoint strong spanning subdigraphs in digraph compositions and digraph products were given. After that, Bang-Jensen, Gutin and Yeo \cite{Bang-Jensen-Gutin-Yeo} characterized all strong semicomplete compositions with a pair of arc-disjoint strong spanning subdigraphs. Recently, Gutin and Sun \cite{Gutin-Sun} studied the existence of a pair (it is called ``good pair'') of arc-disjoint in- and out-branchings rooted at the same vertex in this class of digraphs. Especially, they characterized semicomplete compositions with a good pair rooted at the same vertex, which generalizes the corresponding characterization by Bang-Jensen and Huang \cite{Bang-Jensen-Huang1} for quasi-transitive digraphs. In \cite{Sun}, Sun studied the topic of $k$-kings on compositions of digraphs.
When $T$ is arbitrary, digraph compositions with a $k$-king and digraph compositions all of whose vertices are $k$-kings were characterized.
When $T$ is semicomplete, he discussed the existence of 3-kings and 4-kings, and the adjacency between a 3-king and a non-king.

%we deduce that every strong semicomplete composition has a 3-king. We also show that $Q$ has more 3-kings when $T$ has no vertex of in-degree zero. The adjacency between a 3-king and a non-king in a semicomplete composition are also studied. We then give a sufficient condition under which a strong semicomplete composition can be established, that is, there exists a strong semicomplete composition $Q'$ which contains $Q$ as an induced subdigraph such that the set of all 3-kings of $Q'$ is precisely $V(Q)$. Finally, we discuss the minimum number of 4-kings in a strong semicomplete composition.

In this paper, we will continue the research of digraph compositions and focus on the structural properties of semicomplete compositions. Our results show that this class of digraphs shares some nice structural properties of other classes of digraphs, like quasi-transitive digraphs.

In Section~\ref{sec:connectivity}, connectivity of semicomplete compositions $Q=T[H_1, \dots, H_t]$ will be studied. Let $\mathcal{T}_1$ be the set of all semicomplete digraphs $T$ satisfying the following: there is some vertex $u$ such that $uv, vu\in A(T)$ for any $v\in V(T)\setminus \{u\}$. We will first investigate connectivity of a strong semicomplete composition $Q$ with $T\not \in \mathcal{T}_1$ (Theorems~\ref{thmd} and~\ref{thmi}), and examples show that Theorems~\ref{thmd} and~\ref{thmi} may not hold when $T\in \mathcal{T}_1$. Especially, the structure of a minimal separator in $Q$ will be given (Theorem~\ref{thmd}).  We then deduce that a strong semicomplete composition $Q$ with at least four vertices has two distinct vertices $v_1, v_2$ such that $Q-v_i$ is strong for $i\in [2]$ (Theorem~\ref{thma}), which extends a similar result of quasi-transitive digraphs by Bang-Jensen and Huang \cite{Bang-Jensen-Huang1}.

In Section~\ref{sec:paths}, we will focus on the structure of paths in semicomplete compositions. Semicomplete compositions which have a Hamiltonian path will be characterized (Theorem~\ref{thmg}). In addition, results on $k$-path subdigraphs of semicomplete compositions will be given (Theorems~\ref{thmb} and~\ref{thmb1}), one of them extends a similar result by Bang-Jensen, Nielsen and Yeo \cite{Bang-Jensen-Nielsen-Yeo}.

In Section~\ref{sec:cycles}, we will first give two characterizations (Theorems~\ref{thmg1} and~\ref{thmg2}) of Hamiltonian semicomplete compositions which extend the characterizations of Hamiltonian quasi-transitive digraphs by Bang-Jensen \& Huang \cite{Bang-Jensen-Huang1}, and Gutin \cite{Gutin4}, respectively.
Bang-Jensen and Huang \cite{Bang-Jensen-Huang1} used the similarities between extended semicomplete digraphs and quasi-transitive digraphs to derive results on pancyclic and vertex-pancyclic Hamiltonian quasi-transitive digraphs. We will extend their result on the pancyclicity of quasi-transitive digraphs to semicompete compositions by characterizing all pancyclic semicomplete compositions (Theorem~\ref{thme}).

In Section~\ref{sec:spanning}, we will turn attention to strong spanning subdigraphs in semicomplete compositions. Strong semicomplete compositions which contain a spanning strong semicomplete composition without a 2-cycle will be characterized (Theorem~\ref{thma1}). This is a similar result to that of strong quasi-transitive digraphs by Bang-Jensen and Huang \cite{Bang-Jensen-Huang1}. The {\em smallest strong spanning subdigraph} of a digraph $D$ is defined to be the strong spanning subdigraph of $D$ with smallest arcs. The goal of {\sc Smallest Strong Spanning Subdigraph} (SSSS) problem is finding a smallest strong spanning subdigraph of a given strong digraph $D$. This problem generalizes the Hamiltonian cycle problem (and therefore is NP-hard) and has been considered in the literature, see e.g. \cite{Aho-Garey-Ullman, Bang-Jensen-Huang-Yeo, Khuller-Raghavachari-Young}. Especially, Bang-Jensen, Huang and Yeo \cite{Bang-Jensen-Huang-Yeo} deduced that the smallest strong spanning subdigraph of a strong quasi-transitive digraph $D$ has precisely $n+\epsilon(D)$ arcs, where $\epsilon(D)$ is defined in Section~\ref{sec:spanning}. We will study the smallest strong spanning subdigraph of strong semicomplete compositions by extending this result to strong semicomplete compositions (Theorem~\ref{thmsss}).

In Section~\ref{sec:acyclic}, we will give a sufficient condition to guarantee the existence of an acyclic spanning subdigraph $R$ with a source $x$ and a sink $y$ such that for each vertex $z$ of $Q$, $R$ contains an $x-z$ path and a $z-y$ path (Theorem~\ref{thmf}). Similar results for tournaments and quasi-transitive digraphs were given by Thomassen \cite{Thomassen2}, and Bang-Jensen and Huang \cite{Bang-Jensen-Huang1}, respectively.

\section{Connectivity}\label{sec:connectivity}

By the definition of semicomplete compositions, we directly have the following observation.

\begin{obs}\label{obs1}\cite{Sun}
Let $Q=T[H_1, \dots, H_t]$ be a digraph composition, $Q$ is strong if and only if $T$ is strong.
\end{obs}

Recall that $\mathcal{T}_1$ is the set of all semicomplete digraphs $T$ satisfying the following: there exists a vertex $u$ such that $uv, vu\in A(T)$ for any $v\in V(T)\setminus \{u\}$. The following result concerns the structure of a minimal separator in a strong semicomplete composition $Q=T[H_1, \dots, H_t]$ when $T\not \in \mathcal{T}_1$.

\begin{thm}\label{thmd}
Let $Q=T[H_1, \dots, H_t]$ be a strong semicomplete composition with $T\not \in \mathcal{T}_1$. %Suppose $S$ is a minimal separator of $Q$.
Then every minimal separator $S$ of $Q$ induces in $\overline{U(Q)}$ a subgraph which consists of some connected components of $\overline{U(Q)}$. Moreover, each vertex $s\in S$ is adjacent to every vertex of $Q-S$.
\end{thm}
\begin{pf}
%Since $T\not \cong \overleftrightarrow{K_n}$, we have that each minimal separator of $T$ has size at most $t-2$.
By Observation~\ref{obs1}, $T$ is strong. Furthermore, we have $t\geq 3$ since $T\not \in \mathcal{T}_1$. Let $S$ be a minimal separator of $Q$, we have the following claim:

{\bf Claim.} $Q-S$ contains vertices from at least two distinct $H_i$.

{\bf Proof of the claim.} Suppose that $Q-S\subseteq V(H_i)$ for some $i$, say $i=1$. Observe that now we have $\bigcup_{i=2}^t{V(H_i)}\subseteq S$. Since $T\not\in \mathcal{T}_1$, there exists some vertex, say $u_2$, such that at least one of $u_1u_2$ and $u_2u_1$ does not belong to $A(T)$. By Observation~\ref{obs1}, $H_1\cup H_2=Q-(\bigcup_{i=3}^t{V(H_i)})$ is not strong, which means that $\bigcup_{i=3}^t{V(H_i)}\subsetneqq S$ is a separator of $Q$, a contradiction. Hence, the claim holds.\qed%$Q-S$ contains vertices from at least two distinct $H_i$.

%then $S_T=V(T)\setminus \{v_1\}$ is a minimal separator of $T$. Otherwise, some proper subset, say $S'_T$, of $S_T$ is a separator of $T$. Let $S'$ be the union of all $V(H_i)$ such that $u_i \in S'_T$. By Observation~\ref{obs1}, $Q-S'$ is not strong, so $S'$ is a separator of $Q$ and is a proper subset of $S$, which produces a contradiction. So $S_T$ is a minimal separator of $T$, but now this produces a contradiction with the fact that each minimal separator of $T$ has size at most $t-2$. Hence, $Q-S$ contains vertices from at least two distinct $H_i$ and the claim holds.

Let $H\subseteq V(Q)$ induce a connected component of $\overline{U(Q)}$. Suppose that $H-S\neq \emptyset$ and $H\cap S\neq \emptyset$. Let $Q'=Q-(S-H)$. By the minimality of $S$, the definition of semicomplete compositions and the above claim, the subdigraph $Q'$ is a strong semicomplete composition with $Q'=T'[H'_1, \dots, H'_{t'}]$, where $T'$ is a strong subgraph of $T$ (by Observation~\ref{obs1}) with $|T'|=t'$ and each $H'_i$ is a subdigraph of some $H_j$. Let $Q''=Q'-(H\cap S)=Q-S$. Observe that $Q''$ is also a semicomplete composition. Furthermore, since $H-S\neq \emptyset$, $Q''$ must be a composition of the form $T'[H''_1, \dots, H''_{t'}]$ where each $H''_i$ is a subdigraph of some $H_j$, and hence is also strong by Observation~\ref{obs1}, contradicting that $S$ is a separator. Therefore, either $H\cap S=\emptyset$ or $H\subseteq S$, which means that $S$ induces in $\overline{U(Q)}$ a subgraph which consists of some connected components of $\overline{U(Q)}$.

Let $s\in S$. For any vertex $x\in Q-S$, if $x$ is not in the same $H_i$ as $s$, then $x$ and $s$ are adjacent by the definition of a semicomplete composition. Otherwise, $x$ and $s$ belong to the same $H_i$. If $x$ and $s$ are nonadjacent in $Q$, then they must be in the same component in $\overline{U(Q)}$, therefore $x\in S$ by the above argument, this produces a contradiction.
\end{pf}

Note that the above result may not hold when $T\in \mathcal{T}_1$. Consider the following example: let $T$ be a semicomplete digraph such that for each $1\leq j\leq t$, $u_1u_j, u_ju_1\in A(T)$ and $H_i$ is a cycle of length four: $u_{i,1}, u_{i,2}, u_{i,3}, u_{i,4}, u_{i,1}$. Observe that $S=\{u_{1,1}\}\cup (\bigcup_{i=2}^t{V(H_i)})$ is a minimal separator of $Q$, and $\overline{U(Q)}$ consists of $2t$ components (each component is an edge). Clearly, the edge $\{u_{1,1}, u_{1,3}\}$ is one such component, however, $u_{1,1}\in S$ and $u_{1,3}\not\in S$.

We continue to study the connectivity of semicomplete compositions as follows.

\begin{thm}\label{thmi}
Let $Q=T[H_1, \dots, H_t]$ be a $k$-strong-connected semicomplete composition with $T\not \in \mathcal{T}_1$. If $V(H_i)$ induces a connected component of $\overline{U(Q)}$ for some $i\in [t]$, then deleting all arcs in $H_i$ results in a $k$-strong-connected semicomplete composition.
\end{thm}
\begin{pf} Without loss of generality, assume that $V(H_1)$ induces a connected component of $\overline{U(Q)}$. The result clearly holds for the case that $H_1$ has no arcs, so in the following we assume that $H_1$ has at least one arc. Let $Q'$ be the subdigraph of $Q$ by deleting all arcs in $H_1$. Let $S$ be a minimum separator of $Q'$, we will show that $Q-S$ is not strong and then the result holds.

If $Q'-S$ contains vertices from at least two $H_i$, say $H_{i_1}, \dots, H_{i_{t'}}$, then $Q'-S$ is a semicomplete composition of the form $T'[H'_{i_1}, \dots, H'_{i_{t'}}]$ where $t'=|T'|\geq 2$ and each $H'_{i_j}$ is a subdigraph of $H_{i_j}$ for $j\in [t']$. Furthermore, $Q-S$ is a semicomplete composition of the form $T''[H''_{i_1}, \dots, H''_{i_{t''}}]$ where $t''=|T''|\geq 2$ and each $H''_{i_j}$ is a subdigraph of $H_{i_j}$ for $j\in [t'']$. Clearly, $T''=T'$. By Observation~\ref{obs1}, $T'$ is not strong and so $Q-S$ is also not strong.

Otherwise, $Q'-S$ contains vertices from exactly one $H_i$. We need the following claim.

{\bf Claim.} $i\neq 1$.

{\bf Proof of the claim.} If $i=1$, then $Q'-S=Q'[V(H_1)]$ by the minimality of $S$ and the fact that $Q'[V(H_1)]$ has no arcs. Assume that there exists some $j\in [t]\setminus \{1\}$ such that at most one of $u_1u_j, u_ju_1$ belongs to $A(T)$. Let $Q''=Q'-(S\setminus \{u\})$ where $u\in V(H_j)$. Observe that $Q''$ is not strong and so $S\setminus \{u\}$ is a separator which is a proper subset of $S$, a contradiction. Hence, for each $j\in [t]\setminus \{1\}$, both $u_1u_j$ and $u_ju_1$ belong to $A(T)$, but now we have $T\in \mathcal{T}_1$, this also produces a contradiction.\qed

By the above claim, we have $V(H_1) \subseteq S$, and then $Q-S=Q'-S$ (is not strong). This completes the proof.
\end{pf}

Note that Theorem~\ref{thmi} may not hold if $T\in \mathcal{T}_1$. Consider the following example: let $T$ be a 2-cycle, and $H_i$ be a $2k$-cycle for each $1\leq i\leq 2$, where $k\geq 3$. Clearly, $T\in \mathcal{T}_1$. Observe that each $V(H_i)$ induces a connected component in $\overline{U(Q)}$ for $1\leq i\leq 2$. Indeed, each such component consists of two vertex-disjoint cliques of size $k$ and edges between them. Let $Q'$ be the subdigraph of $Q$ by deleting all arcs in $H_1$. It can be checked that $V(H_2)$ is a separator of size $k$ in $Q'$. However, there is no separator of size $k$ in $Q'$. Hence, $Q$ is $(k+1)$-strong-connected, but $Q'$ is not $(k+1)$-strong-connected.

The following result on strong quasi-transitive digraphs was given by Bang-Jensen and Huang.

\begin{thm}\label{pro02}\cite{Bang-Jensen-Huang1}
Every strong quasi-transitive digraph $D$ with at least four vertices has two distinct vertices $v_1, v_2$ such that $D-v_i$ is strong for $i\in [2].$
\end{thm}

We extend Theorem~\ref{pro02} to strong semicomplete compositions.

\begin{thm}\label{thma}
Every strong semicomplete composition $Q=T[H_1, \dots, H_t]$ with at least four vertices has two distinct vertices $v_1, v_2$ such that $Q-v_i$ is strong for $i\in [2].$
\end{thm}
\begin{pf} Let $Q=T[H_1, \dots, H_t]$ be a strong semicomplete composition.
For the case that $2\leq t\leq 3$, since $Q$ has at least four vertices, we clearly have that there exists one $H_i$ with at least two vertices, say $v_1$ and $v_2$. For each $i\in [2]$, the digraph $Q-v_i$ is a semicomplete composition of the form $T'[H'_1, \dots, H'_t]$ where $T'=T$. Hence, by Observation~\ref{obs1}, $Q-v_i$ is strong.

Now we consider the case that $t\geq 4$. If each $H_i$ is trivial, then $Q=T$ is a semicomplete digraph and the result holds by Theorem~\ref{pro02}. Otherwise, there exists one $H_i$ with at least two vertices, then with a similar argument to that of the above paragraph, the result holds.
\end{pf}

\section{Paths}\label{sec:paths}

\subsection{Hamiltonian paths}

Gutin proved the following result on Hamiltonian paths and cycles of extended semicomplete digraphs.

\begin{thm}\label{Gutin-Hamiltonian}\cite{Bang-Jensen-Gutin-Huang}
An extended semicomplete digraph has a Hamiltonian path (resp. cycle) if and only if it has a 1-path-cycle factor (resp. it is strong and has a cycle factor).
\end{thm}

By the above theorem, we can characterize semicomplete compositions which have a Hamiltonian path.

\begin{thm}\label{thmg}
Let $Q=T[H_1, \dots, H_t]$ be a semicomplete composition. Then
$Q$ has a Hamiltonian path if and only if it has a 1-path-cycle
factor $\mathcal{F}=P\cup C_1\cup \dots \cup C_k~(k\geq 0)$ such that neither $V(P)$ nor $V(C_i)$ is completely contained in a connected component of $\overline{U(Q)}$.
\end{thm}
\begin{pf}
%We just prove the assertion $(a)$ since the argument of $(b)$ is similar.
Suppose that $Q$ has a 1-path-cycle factor $\mathcal{F}=P\cup C_1\cup \dots \cup C_k~(k\geq 0)$, such that neither $V(P)$ nor $V(C_i)$ is completely contained in a connected component of $\overline{U(Q)}$. Let $V_1, V_2, \dots, V_r$ be the partition of $Q$ such that each $V_i$ induces a connected component of $\overline{U(Q)}$. Observe that $r\geq 2$ by the definition of a semicomplete composition. We obtain a digraph $Q'$ from $Q$ as follows: Shrink each maximal subpath of $P$ and $C_i$ (for each $1\leq i\leq k$) which is completely contained in some $Q[V_j]$, that is, each such subpath is a vertex of $Q'$; delete all arcs which are still in $Q[V_j]$. It can be checked that $Q'$ is an extended semicomplete digraph and has a 1-path-cycle factor. By Theorem~\ref{Gutin-Hamiltonian}, $Q'$ has a Hamiltonian path which can be transformed to a Hamiltonian path of $Q$. The other direction is clearly true.
\end{pf}

\subsection{Path subdigraphs}

The following structural characterization of longest cycles in extended semicomplete digraphs was extensively used.

\begin{thm}\label{lem05}\cite{Bang-Jensen-Huang-Yeo}
Let $Q=T[\overline{K}_{n_1}, \dots, \overline{K}_{n_t}]$ be a strong extended semicomplete digraph. For
$i\in [t]$, let $m_i$ denote the maximum number of vertices from $\overline{K}_{n_i}$ which can be covered by a cycle subdigraph of $Q$. Then every longest cycle of $Q$ contains precisely $m_i$ vertices from $\overline{K}_{n_i}$ for each $i\in [t]$.
\end{thm}

Bang-Jensen et al. proved a similar result for $k$-path subdigraphs in an extended semicomplete digraph.

\begin{thm}\label{lem03}\cite{Bang-Jensen-Nielsen-Yeo}
Let $Q=T[\overline{K}_{n_1}, \dots, \overline{K}_{n_t}]$ be an extended semicomplete digraph and let $\ell_{i,k}$ denote the maximum number of vertices of $\overline{K}_{n_i}$ that can be covered by a $k$-path subdigraph in $Q$. Then every maximum $k$-path subdigraph in $Q$ covers exactly $\ell_{i,k}$ vertices of $\overline{K}_{n_i}$ for each $i\in [t]$.
\end{thm}

By Theorem~\ref{lem03}, Bang-Jensen et al. got the following result on a strong semicomplete composition when each $H_i$ is either a single vertex or a non-strong quasi-transitive digraph.

\begin{thm}\label{lem04}\cite{Bang-Jensen-Nielsen-Yeo}
Let $Q=T[H_1, \dots, H_t]$ be a strong semicomplete composition with order $n$, where each $H_i$ is either a single vertex or a non-strong quasi-transitive digraph. For every $k\in [n]$ and $i\in [t]$, there exists an integer $n_{i,k}$ such that every maximum $k$-path subdigraph $\mathcal{F}$ of $Q$ satisfies $|V(H_i) \cap V(\mathcal{F})|=n_{i,k}$ and no $k$-path subdigraph of $Q$ contains more than $n_{i,k}$ vertices of $H_i$.
\end{thm}

We extend Theorem~\ref{lem04} to all strong semicomplete compositions.%in the following result which also means that Theorem~\ref{lem03} holds for a general strong semicomplete composition.

\begin{thm}\label{thmb}
Let $Q=T[H_1, \dots, H_t]$ be a strong semicomplete composition with order $n$.
For every $k\in [n]$ and $i\in [t]$, there exists an integer $n_{i,k}$ such that every maximum $k$-path subdigraph $\mathcal{F}$ of $Q$ satisfies $|V(H_i) \cap V(\mathcal{F})|=n_{i,k}$ and no $k$-path subdigraph of $Q$ contains more than $n_{i,k}$ vertices of $H_i$.
\end{thm}
\begin{pf}
Let $Q'= T[\overline{K}_{n_1}, \dots, \overline{K}_{n_t}]$ be the corresponding extended semicomplete digraphs of $Q$. By Observation~\ref{obs1}, $T$ is strong and so $Q'$ is a strong extended semicomplete digraph. For each $i\in [t]$, let $\ell_{i,k}$ be the maximum number of vertices of $\overline{K}_{n_i}$ that can be covered by a $k$-path subdigraph in $Q'$; moreover, let $n_{i,k}$ be the maximum number of vertices in an $\ell_{i,k}$-path subdigraph in $H_i$. By Theorem~\ref{lem03}, there exists a $k$-path subdigraph $\mathcal{F}'$ in $Q'$ containing exactly $\ell_{i,k}$ vertices from $\overline{K}_{n_i}$ for each $i\in [t]$. Now we obtain a $k$-path subdigraph $\mathcal{F}$ in $Q$ from $\mathcal{F}'$ by replacing each of the $\ell_{i,k}$ vertices of $\overline{K}_{n_i}$ by a path from an $\ell_{i,k}$-path subdigraph of $H_i$ which covers $n_{i,k}$ vertices.

Furthermore, we claim that no $k$-path subdigraph in $Q$ can include more than $n_{i,k}$ vertices from $H_i$. Indeed, if there is such a $k$-path subdigraph, say $\mathcal{F}_0$, in $Q$, then $\mathcal{F}_0$ must visit $H_i$ $\ell'_{i,k}~(> \ell_{i,k})$ times, that is, there are $\ell'_{i,k}$ disjoint paths in $H_i$ such that each such path is a subpath of an element of $\mathcal{F}_0$. Now we can obtain a $k$-path subdigraph, say $\mathcal{F}'_0$, from $\mathcal{F}_0$ by shrinking all subpaths of $\mathcal{F}_0$ contained in $H_j$ for each $j\in [t]$. Since all vertices in the same $H_i$ are similar, we could see $\mathcal{F}'_0$ as a $k$-path subdigraph of $Q'$ which visits $\overline{K}_{n_i}$ $\ell'_{i,k}~(> \ell_{i,k})$ times, this produces a contradiction. Therefore, the conclusion holds.%every maximum $k$-path subdigraph $\mathcal{F}$ of $Q$ satisfies $|V(H_i) \cap V(\mathcal{F})|=n_{i,k}$ for $i\in [t]$.
\end{pf}

The argument for Theorem~\ref{thmb} also means that Theorem~\ref{lem03} holds for a strong semicomplete composition. Note that the number $\ell_{i, k}$ in Theorem~\ref{thmb1} is exactly $n_{i, k}$ of Theorem~\ref{thmb}.

\begin{thm}\label{thmb1}
Let $Q=T[H_1, \dots, H_t]$ be a strong semicomplete composition and let $\ell_{i,k}$ denote the maximum number of vertices of $H_i$ that can be covered by a $k$-path subdigraph in $Q$. Then every maximum $k$-path subdigraph in $Q$ covers exactly $\ell_{i,k}$ vertices of $H_i$ for each $i\in [t]$.
\end{thm}

\section{Cycles}\label{sec:cycles}

\subsection{Hamiltonian cycles}

With a similar argument to that of Theorem~\ref{thmg}, we provide the first characterization of Hamiltonian semicomplete compositions which extends a characterization of Hamiltonian quasi-transitive digraphs given by Bang-Jensen and Huang \cite{Bang-Jensen-Huang1}.

\begin{thm}\label{thmg1}
Let $Q=T[H_1, \dots, H_t]$ be a semicomplete composition. Then $Q$ has a Hamiltonian cycle if and only if it is strong and contains a cycle factor
$\mathcal{F}=C_1\cup \dots \cup C_k$, such that no $V(C_i)$ is completely contained in a connected component of $\overline{U(Q)}$.
\end{thm}

%The following characterization of Hamiltonian quasi-transitive digraphs is given implicitly in \cite{Gutin4} by Gutin:

%\begin{thm}\label{Gutin-Hamil}\cite{Gutin4}
%Let $D$ be a strong quasi-transitive digraph with canonical decomposition $D= S[Q_1, \dots, Q_s]$. Let $n_i$ be the order of the digraph $Q_i$ for $i\in [s]$. Then $D$ is Hamiltonian if and only if the extended semicomplete digraph $D'=S[\overline{K}_{n_1}, \dots, \overline{K}_{n_s}]$ has a cycle subdigraph which covers at least $pc(Q_i)$ vertices of $\overline{K}_{n_i}$ for every $i \in [s]$.
%\end{thm}

Gutin obtained the following result on long cycles in strong extended semicomplete digraphs.

\begin{thm}\label{Gutin-longcycle}\cite{Gutin2}
Let $D$ be a strong extended semicomplete digraph and let $\mathcal{F}$ be a cycle subdigraph of $D$. Then $D$ has a cycle $C$ which contains all vertices of $\mathcal{F}$. %The cycle $C$ can be found in time $O(|V(D)|^3)$.
In particular, if $V(\mathcal{F})$ is maximum, then $V(C)=V(\mathcal{F})$ and $C$ is a longest cycle of $D$.
\end{thm}

We now extend a characterization of Hamiltonian quasi-transitive digraphs which is given by Gutin \cite{Gutin4} to strong semicomplete compositions, and provide the second characterization of Hamiltonian semicomplete compositions.

\begin{thm}\label{thmg2}
Let $Q=T[H_1, \dots, H_t]$ be a strong semicomplete composition. Then $Q$ has a Hamiltonian cycle if and only if the extended semicomplete digraph $Q'=T[\overline{K}_{n_1}, \dots, \overline{K}_{n_t}]$ has a cycle subdigraph which covers at least $pc(H_i)$ vertices of $\overline{K}_{n_i}$ for every $i \in [t]$.
\end{thm}
\begin{pf}
Suppose that $Q$ has a Hamiltonian cycle $C$. For each $i\in [t]$, the subdigraph $V(H_i)\cap C$ is a $k_i$-path factor, denoted by $\mathcal{F}_i$, of $H_i$, where $k_i$ is a positive integer. Observe that $k_i\geq pc(H_i)$ by the definition of the path covering number. We can transform $C$ into a cycle of $Q'$ which covers at least $pc(H_i)$ vertices of $\overline{K}_{n_i}$ for every $i \in [t]$, by the following operations: For each $i\in [t]$, delete the arcs between end-vertices of all paths in $\mathcal{F}_i$ except for the paths themselves and then shrink all paths in $\mathcal{F}_i$.

Now we suppose that $Q'$ has a cycle subdigraph $\mathcal{F}'$ containing $k_i(\geq pc(H_i))$ vertices of $\overline{K}_{n_i}$ for every $i \in [t]$. By Observation~\ref{obs1} and the assumption that $Q$ is strong, $Q'$ is strong and hence has a cycle $C'$ such that $V(C')=V(\mathcal{F}')$ by Theorem~\ref{Gutin-longcycle}. Observe that $H_i$ has a $k_i$-path factor $\mathcal{F}_i$. We can obtain a Hamiltonian cycle in $Q$ by replacing the $k_i$ vertices of $\overline{K}_{n_i}$ in $C'$ with the paths of $\mathcal{F}_i$ for each $i\in [t]$.
\end{pf}

\subsection{Pancyclicity}%\label{sec:pancyclicity}

Gutin characterized extended semicomplete digraphs which are pancyclic. %and vertex-pancyclic extended semicomplete digraphs.
\begin{thm}\label{thm05}\cite{Gutin}
Let $D$ be a Hamiltonian extended semicomplete digraph on $n\geq 4$ vertices such that $\overline{U(D)}$ has $k\geq 3$ connected components. Then $D$ is pancyclic if and only if $D$ is not triangular with a partition $\{V_0, V_1, V_2\}$, two of which induce independent sets, such that either $|V_0|=|V_1|=|V_2|$ or no $D[V_i]$ $(i= 0, 1, 2)$ contains a path of length 2.
%the following assertions hold:
%\begin{description}
%\item[(a)] $D$ is pancyclic if and only if $D$ is not triangular with a partition $\{V_0, V_1, V_2\}$, two of which induce independent sets, such that either $|V_0|=|V_1|=|V_2|$ or no $D[V_i]$ $(i= 0, 1, 2)$ contains a path of length 2.
%\item[(b)] $D$ is vertex-pancyclic if and only if it is pancyclic and either $k>3$ or $k = 3$ and $D$ contains two cycles $C, C'$ of length 2 such that $C \cup C'$ has vertices in each of the three connected components of $\overline{U(D)}$.
%\end{description}
\end{thm}

The next two lemmas by Bang-Jensen and Huang concern cycles in triangular digraphs.

\begin{lem}\label{lem01}\cite{Bang-Jensen-Huang1}
Suppose that $D$ is a Hamiltonian triangular digraph with a partition $\{V_0, V_1, V_2\}$. If $D[V_1]$ contains an arc $xy$ and $D[V_2]$ contains an arc $uv$, then every vertex of $V_0\cup \{x, y, u, v\}$ is on cycles of lengths $3, 4, \dots, |V(D)|$.
\end{lem}

\begin{lem}\label{lem02}\cite{Bang-Jensen-Huang1}
Suppose that $D$ is a triangular digraph with a partition $\{V_0, V_1, V_2\}$ and has a Hamiltonian cycle $C$. If $D[V_0]$ contains an arc of $C$ and a path $P$ of length 2, then every vertex of $V_1\cup V_2\cup V(P)$ is on cycles of lengths $3, 4, \dots, |V(D)|$.
\end{lem}

By Theorem~\ref{thm05}, Lemmas~\ref{lem01} and \ref{lem02}, we get the following result on pancyclicity of semicomplete compositions, and this extends a similar result for quasi-transitive digraphs by Bang-Jensen and Huang \cite{Bang-Jensen-Huang1}.

\begin{thm}\label{thme}
Let $Q=T[H_1, \dots, H_t]$ be a Hamiltonian semicomplete composition on $n\geq 4$ vertices. Then $Q$ is pancyclic if and only if it is not triangular with a partition
$\{V_0, V_1, V_2\}$, two of which induce independent sets, such that either $|V_0|=|V_1|=|V_2|$, or no $Q[V_i]~(i=0, 1, 2)$ contains a path of length 2.
\end{thm}
\begin{pf}
We first prove the necessity of the theorem. Suppose that $Q$ is triangular with a partition $\{V_0, V_1, V_2\}$, two of which induce independent sets. If $|V_0|=|V_1|=|V_2|$, then $Q$ contains no cycle of length $n-1$. If no $Q[V_i]~(i = 0, 1, 2)$ contains a path of length 2, then $Q$ contains no cycle of length 5.

We next prove the sufficiency of the theorem. Let $C$ be a Hamiltonian cycle of $Q$. We obtain an extended semicomplete digraph $Q'$ from $Q$ as follows: For each $i\in [t]$, shrink each subpath of $C$ which is contained in $H_i$, and then delete the remaining arcs (if exist) of $H_i$. Observe that in this process $C$ is transformed to a Hamiltonian cycle $C'$ of $Q'$.

Suppose $Q$ is not pancyclic. It is not difficult to see that $Q'$ is not pancyclic. Now $Q'$ is triangular with a partition $\{V'_0, V'_1, V'_2\}$ by Theorem~\ref{thm05}. For each $i= 0, 1, 2$, let $V_i \subseteq V(Q)$ be obtained from $V'_i$ by substituting back all vertices on shrunk subpaths of $C$. Then $Q$ is triangular with partition $\{V_0, V_1, V_2\}$, moreover, each $Q[V_i]$ is covered by $\ell$ disjoint subpaths of $C$ for some $\ell$.

By Lemma~\ref{lem01}, two elements of $\{V_0, V_1, V_2\}$, say $V_1$ and $V_2$, induce independent sets in $Q$ (Indeed, suppose that $V_1$ and $V_2$ are not independent, then both of them contain at least one arc. By Lemma~\ref{lem01}, every vertex of $V_0\cup \{x, y, u, v\}$ is on cycles of lengths $3, 4, \dots, |V(Q)|$, but this means that $Q$ is pancyclic, a contradiction). If $V_0$ is not independent and $|V_0| > |V_1|$, then $Q[V_0]$ contains an arc of $C$. By Lemma~\ref{lem02}, there is no path of length 2 in $Q[V_0]$. This completes the proof.
\end{pf}

\section{Strong spanning subdigraphs}\label{sec:spanning}

In this section, we concern strong spanning subdigraphs without a 2-cycle and smallest spanning subdigraphs in a semicomplete composition.

\subsection{Strong spanning subdigraphs without a 2-cycle}

The following result on strong semicomplete digraphs can be found in the literature, see e.g. Proposition~2.2.8 of \cite{Bang-Jensen-Havet}.

\begin{thm}\label{pro01}
Every strong semicomplete digraph on $n \geq 3$ vertices contains a strong spanning tournament.
\end{thm}

Theorem~\ref{pro01} means that every strong semicomplete digraph on $n \geq 3$ vertices contains a strong spanning semicomplete digraph without a 2-cycle. Bang-Jensen and Huang obtained a similar result for strong quasi-transitive digraphs.

\begin{thm}\label{pro03}\cite{Bang-Jensen-Huang1}
Every strong quasi-transitive digraph contains a strong spanning quasi-transitive digraph without a 2-cycle.
\end{thm}

We get a similar result to Theorem~\ref{pro03} for strong semicomplete compositions.

\begin{thm}\label{thma1}
A strong semicomplete composition $Q=T[H_1, \dots, H_t]$ contains a strong spanning semicomplete composition without a 2-cycle if
and only if $t\geq 3$.
\end{thm}
\begin{pf} Let $Q=T[H_1, \dots, H_t]$ be a strong semicomplete composition.
If $t=2$, then $T$ is a 2-cycle by Observation~\ref{obs1}. Hence, any spanning strong semicomplete composition without 2-cycle contained in $Q$ must be of the form $Q'=T'[H'_1, H'_2]$, where $T'$ is a spanning subdigraph of $T$ with $T'\neq T$ and each $H'_i$ is a spanning subdigraph of $H_i$ for $1\leq i\leq 2$. However, it is impossible, since now $T'$ is not strong by Observation~\ref{obs1}.

For the case that $t\geq 3$, by Observation~\ref{obs1}, $T$ is strong and therefore contains a strong spanning tournament $T'$ with order $t\geq 3$ by Theorem~\ref{pro01}. Furthermore, by Observation~\ref{obs1}, the semicomplete composition $Q'=T'[H'_1, \dots, H'_t]$ is a strong spanning subdigraph of $Q$ and has no 2-cycle, where $H'_i$ is obtained from $H_i$ by deleting all arcs for each $1\leq i\leq t$. Hence, $Q$ contains a spanning strong semicomplete composition without a 2-cycle if and only if $t\geq 3$.
\end{pf}

\subsection{Smallest strong spanning subdigraph}

An {\em ear decomposition} of a digraph $D$ is a sequence
$\mathcal{P}=(P_0, P_1, P_2, \cdots, P_t)$, where $P_0$ is a cycle
or a vertex and each $P_i$ is a path or a cycle with the following
properties:\\
$(a)$~$P_i$ and $P_j$ are arc-disjoint when $i\neq j$.\\
$(b)$~For each $i=0,1,2,\cdots,t$: let $D_i$ denote the digraph with
vertices $\bigcup_{j=0}^i{V(P_j)}$ and arcs
$\bigcup_{j=0}^i{A(P_j)}$. If $P_i$ is a cycle, then it has
precisely one vertex in common with $V(D_{i-1})$. Otherwise the end
vertices of $P_i$ are distinct vertices of $V(D_{i-1})$ and no other
vertex of $P_i$ belongs to $V(D_{i-1})$.\\
$(c)$~$\bigcup_{j=0}^t{A(P_j)}=A(D)$.

The following result is well-known, see e.g., \cite{Bang-Jensen-Gutin}.
\begin{thm}\label{thmed}
Let $D$ be a digraph with at least two vertices. Then $D$ is strong
if and only if it has an ear decomposition. Furthermore, if $D$ is
strong, every cycle can be used as a starting cycle $P_0$ for an ear
decomposition of $D$.%, and there is a linear-time algorithm to find such
an ear decomposition.
\end{thm}

For a digraph $D$ and a natural number $k$, let $H_k(D)$ denote a digraph obtained from $D$ as follows: add two sets of $k$ new vertices $x_1, x_2, \dots, x_k$, $y_1, y_2, \dots, y_k$; add all possible arcs from $V(D)$ to $x_i$ along with all possible arcs from $y_i$ to $V(D)$, $i\in [k]$; add all arcs of the kind $x_iy_j$, $i, j \in [k]$. Clearly, $H_0(D)= D$. Let $D$ be a strong connected digraph and let $\epsilon(D)$ be the smallest $k \geq 0$ such that $H_k(D)$ is Hamiltonian. Observe that %$\epsilon(D)=pc(D)$ when $\epsilon(D) \geq 1$.
\[
\epsilon(D)=\left\{
   \begin{array}{ll}
     0, &\mbox {$D$ is Hamiltonian;}\\
     pc(D), &\mbox {Otherwise.}
   \end{array}
   \right.
\]

Bang-Jensen, Huang and Yeo obtained the following result on the smallest strong spanning subdigraph of a strong quasi-transitive digraph.

\begin{thm}\label{BHY2003}\cite{Bang-Jensen-Huang-Yeo}
The smallest strong spanning subdigraph of a strong quasi-transitive digraph $D$ has precisely $n+\epsilon(D)$ arcs.
%Let $D$ be a strong quasi-transitive digraph. Then every strong spanning subdigraph of $D$ has at least $n+\epsilon(D)$ arcs. Moreover, the smallest strong spanning subdigraph of $D$ has precisely $n+\epsilon(D)$ arcs.
\end{thm}

In the argument of Theorem~\ref{BHY2003}, the authors used the following two lemmas which will also be used in our theorem below.

\begin{lem}\label{BHY2003-1}\cite{Bang-Jensen-Huang-Yeo}
If $D$ is an acyclic extended semicomplete digraph, then $pc(D)= \max\{|I| \mid I~is~an~independent~set~in~D\}$. Furthermore, starting from $D$, one can obtain a path covering with $pc(D)$ paths by removing the vertices of a longest path $pc(D)$ times.
\end{lem}

\begin{lem}\label{BHY2003-2}\cite{Bang-Jensen-Huang-Yeo}
Let $D$ be a strong extended semicomplete digraph and let $C$ be a longest
cycle in $D$. Then $D-V(C)$ is acyclic.
\end{lem}

We now extend Theorem~\ref{BHY2003} to strong semicomplete compositions.
%We give a similar result for strong semicomplete compositions.

\begin{thm}\label{thmsss}
The smallest strong spanning subdigraph of a strong semicomplete composition $Q=T[H_1, \dots, H_t]$ has precisely $n+\epsilon(Q)$ arcs.
%Let $Q=T[H_1, \dots, H_t]$ be a strong semicomplete composition. Then every strong spanning subdigraph of $Q$ has at least $n+\epsilon(Q)$ arcs. Moreover, the smallest strong spanning subdigraph of $Q$ has precisely $n+\epsilon(Q)$ arcs. %Moreover, the bound can be attained.
\end{thm}
\begin{pf}
We first prove the following claim.%first part that every strong spanning subdigraph of $Q$ has at least $n+\epsilon(Q)$ arcs.

{\bf Claim 1.} Every strong spanning subdigraph of $Q$ has at least $n+\epsilon(Q)$ arcs.

{\bf Proof of Claim 1.} Let $Q'$ be a strong spanning subdigraph of $Q$ with $n+k$ arcs such that no proper subdigraph of $Q'$ is spanning and strong. %(by deleting some arcs if necessary). %\textcolor{red}{Clearly, $Q'$ is of the form $Q'=T[H'_1, \dots, H'_t]$ such that each $H'_i$ is a spanning subdigraph of $H_i$ for each $i\in [t]$}.
By Theorem~\ref{thmed}, $Q'$ can be decomposed into a cycle $P_0$ and $k$ paths or cycles $P_1, \dots, P_k$ (each $P_i$ has length at least two by the minimality assumption on $Q'$) with the following
properties:\\
$(a)$~$P_i$ and $P_j$ are arc-disjoint when $i\neq j$.\\
$(b)$~For each $i=0, 1, 2, \dots, k$, let $Q_i$ denote the digraph with vertices $\bigcup_{j=0}^i{V(P_j)}$ and arcs $\bigcup_{j=0}^i{A(P_j)}$. If $P_i$ is a cycle, then it has precisely one vertex in common with $V(Q_{i-1})$; otherwise, the end vertices of $P_i$ are distinct vertices of $V(Q_{i-1})$ and no other vertex of $P_i$ belongs to $V(Q_{i-1})$.\\
$(c)$~$\bigcup_{j=0}^k{A(P_j)}=A(Q')$. Furthermore, this decomposition can be started with $P_0$ as any cycle in $Q'$. It follows that we may choose $P_0$ such that $V(P_0)\not\subseteq V(H_i)$ for $i\in [t]$.

Now $H_k(D)$ has a cycle factor consisting of $P_0$ and $k$ cycles of the form $C_k=y_iP'_ix_iy_i$, where $i\in [k]$ and $P'_i$ is the path one obtains from $P_i$ by removing the vertices it has in common with $V(Q_{i-1})$. Observe that $H_k(D)$ is also a semicomplete composition and no element of $\{P_0, C_i\mid i\in [k]\}$ is completely contained in a connected component of $\overline{U(H_k(D))}$. Then by Theorem~\ref{thmg1}, $H_k(D)$ has a Hamiltonian cycle and so $\epsilon(Q)\leq k$. Therefore, every strong spanning subdigraph of $Q$ has at least $n+\epsilon(Q)$ arcs. This completes the proof of the claim. \qed%Theorem~\ref{BHY2003} means that the bound is sharp when $Q$ is a strong quasi-transitive digraph.

%Next we prove the second part that the smallest strong spanning subdigraph of $Q$ has precisely $n+\epsilon(Q)$ arcs.

If $Q$ is Hamiltonian, then any Hamiltonian cycle is the smallest strong spanning subdigraph of $Q$ has precisely $n+\epsilon(Q)=n$ arcs. In the following, we assume that $Q$ is not Hamiltonian.

Let $Q_0=T[\overline{K}_{n_1}, \dots, \overline{K}_{n_t}]$ be an extended semicomplete digraph which is obtained from $Q$ by deleting all arcs inside $H_i$ for each $i\in [t]$. By Theorem~\ref{thmg2}, $Q_0$ has no cycle subdigraph which covers at least $pc(H_i)$ vertices of $\overline{K}_{n_i}$ for every $i\in [t]$.
For each $i\in [t]$, let $m_i$ denote the maximum number of vertices from $\overline{K}_{n_i}$ which can be covered by a cycle subdigraph of $Q_0$. By Theorem~\ref{lem05}, every longest cycle $C$ of $Q_0$ contains precisely $m_i$ vertices from $\overline{K}_{n_i}$ for each $i\in [t]$. Observe that $pc(H_i)\geq m_i$ for each $i\in [t]$. Let $k=\max\{pc(H_i)-m_i\mid i\in [t]\}$.

Let $Q'_0=T[\overline{K}_{m'_1}, \dots, \overline{K}_{m'_t}]$ with $V(\overline{K}_{m'_i})=\{x_{i, j_i}\mid j_i\in [m'_i]\}$ be an extended semicomplete subdigraph of $Q_0$, where $m'_i=\max\{pc(H_i), m_i\}$ for $i\in [t]$. Since vertices inside an independent set of an extended semicomplete subdigraph are similar, we may think of $C$ as a longest cycle in $Q'_0$, that is, $C$ contains precisely $m_i$ vertices from $\overline{K}_{m'_i}$ for $i\in [t]$. By Lemmas~\ref{BHY2003-1} and~\ref{BHY2003-2}, $Q'_0-V(C)$ can be covered by $k$ paths, say $P'_1, \dots, P'_k$. Since $Q'_0-V(C)$ is acyclic, we may assume (by Lemma~\ref{BHY2003-1}) that $P'_i$ starts at a vertex $x$ and ends at a vertex $y$ such that $x$ has in-degree zero and $y$ has out-degree zero in $Q'_0-V(C)$. Now in $Q'_0$ there exists two vertices $z, z'$ such that $zx, yz'\in A(Q'_0)$, and we can add $P'_1$ to $C$ by adding the arcs $zx, yz'$. Then we obtain a strong spanning subdigraph $Q''_0$ of $Q'_0$ with $|V(Q'_0)|+k$ arcs by induction on $k$ that adding $P'_2, \dots, P'_k$ one by one, using two new arcs each time.

Recall that $m'_i=\max\{pc(H_i), m_i\}$, we have $m'_i\geq pc(H_i)$. By the definition of the path covering number, for each $i\in [t]$, $H_i$ has a set of $m'_i$ disjoint paths, say $P_{i, j_i}~(j_i\in [m'_i])$, which cover all vertices of $H_i$. Now we obtain a strong spanning subdigraph $Q'$ of $Q$ by replacing $x_{i, j_i}$ in $Q''_0$ by the path $P_{i, j_i}$ for each $i\in [t], j_i\in [m'_i]$. Furthermore, we have
\begin{eqnarray*}
|A(Q')|
&=&\sum_{i=1}^t(|V(H_i)|-m'_i)+(|V(Q'_0)|+k) \\[6pt]
&=& (n-|V(Q'_0)|)+(|V(Q'_0)|+k)\\[6pt]
&=&n+k.
\end{eqnarray*}
%$$|A(Q')|=\sum_{i=1}^t(|V(H_i)|-m'_i)+(|V(Q'_0)|+k)=(n-|V(Q'_0)|)+(|V(Q'_0)|+k)=n+k.$$

To finish our argument, we still need the following claim.

{\bf Claim 2.} $\epsilon(Q)\geq k$.

{\bf Proof of Claim 2.} Suppose that $\epsilon(Q)< k$. By the definition of $\epsilon(Q)$, $H_{\epsilon(Q)}(Q)=T'[H_1, \dots, H_t, \overline{K}_{\epsilon(Q)}, \overline{K}_{\epsilon(Q)}]$ has a Hamiltonian cycle $C$, where $T'$ is obtained from $T$ by adding two new vertices $u_{t+1}$, $u_{t+2}$ such that $u_{t+1}u_{t+2}$ is an arc and $u_{t+1}$ is dominated by all vertices of $T$ and $u_{t+2}$ dominates all vertices of $T$.
Let $H'_{\epsilon(Q)}(Q)$ be a subdigraph of $H_{\epsilon(Q)}(Q)$ by the following operations: shrink each subpath of $C$ which lies entirely inside some $H_i$ (the resulting cycle is denoted by $C'$), and delete all remaining arcs (if exist) inside each $H_i$. Observe that $H'_{\epsilon(Q)}(Q)$ must be a semicomplete composition such that $H'_{\epsilon(Q)}(Q)=T'[\overline{K}_{a_1}, \dots, \overline{K}_{a_t}, \overline{K}_{\epsilon(Q)}, \overline{K}_{\epsilon(Q)}]$, moreover, $n_i\geq a_i\geq pc(H_i)$ for each $i\in [t]$ and $C'$ is a Hamiltonian cycle in $H'_{\epsilon(Q)}(Q)$.

Let $\{x_i\mid i\in [\epsilon(Q)]\}$ and $\{y_i\mid i\in [\epsilon(Q)]\}$ be the vertices of $H_{\epsilon(Q)}(Q)$ corresponding to $u_{t+1}$ and $u_{t+2}$, respectively. Now by removing all vertices of $\{x_i, y_i\mid i\in [\epsilon(Q)]\}$, we can obtain a set of $\epsilon(Q)$ disjoint paths, $P_1, \dots, P_{\epsilon(Q)}$, covering all vertices in $H''_{\epsilon(Q)}(Q)=T[\overline{K}_{a_1}, \dots, \overline{K}_{a_t}]$ since all arcs leaving $x_i$ go to $\{y_i\mid i\in [\epsilon(Q)]\}$. Recall that $Q_0=T[\overline{K}_{n_1}, \dots, \overline{K}_{n_t}]$ is an extended semicomplete digraph which is obtained from $Q$ by deleting all arcs inside each $H_i$ for $i\in [t]$. Since all vertices in the same independent set of an extended semicomplete digraph are similar, we can assume that $P_i$ is a path in $Q_0$ for each $i\in [\epsilon(Q)]$. Then $\mathcal{F}=P_1\cup \dots\cup P_{\epsilon(Q)}$ is an $\epsilon(Q)$-path subdigraph of $Q_0$ which covers $a_i$ vertices of $\overline{K}_{n_i}$ in $Q_0$ for each $i\in [\epsilon(Q)]$.

Let $i_0$ be chosen such that $k=pc(H_{i_0})-m_{i_0}$. Recall that $\epsilon(Q)< k$ and $a_i\geq pc(H_i)$ for each $i\in [t]$. We directly have $a_{i_0}\geq pc(H_{i_0})=k+m_{i_0}>\epsilon(Q)+m_{i_0}$ and so $a_{i_0}>\epsilon(Q)$, which means that some path $P_j$ in $\mathcal{F}$ (there maybe more than one such path) contains at least two vertices of $\overline{K}_{n_{i_0}}$ in $Q_0$. Let $P_j=z_1, \dots, z_p$, and $z_a, z_b$ are similar where $a<b$. Observe that $C_j=z_{a+1}, \dots, z_{b-1}, z_b, z_{a+1}$ is a cycle and $z_az_{b+1}\in A(Q_0)$ when $b<p$. Moreover, we can replace $P_j$ by the cycle $C_j$ and a path $P'_j=P_j[z_1, z_a]P_j[z_{b+1}z_p]$ since $V(P_j)=V(C_j)\cup P'_j$. Now we obtain a cycle subgraph of $Q_0$ by continuing the above process (replacing paths in $\mathcal{F}$ which contains at least two vertices of $\overline{K}_{n_{i_0}}$ in $Q_0$ by a cycle and a path) until every path in $\mathcal{F}$ contains at most one vertex from $\overline{K}_{n_{i_0}}$ in $Q_0$. However, this cycle subdigraph covers at least $a_{i_0}-\epsilon(Q)> m_{i_0}$ vertices of $Q_0$, which contradicts the definition of $m_{i_0}$. Therefore, we have $\epsilon(Q)\geq k$.\qed

By Claims 1 and 2, $Q'$ is the smallest strong spanning subdigraph of $Q$ and has precisely $n+\epsilon(Q)$ arcs. This completes the proof.
\end{pf}

\section{Acyclic spanning subdigraphs}\label{sec:acyclic}

We now turn attention to the existence of prescribed acyclic spanning subdigraphs in a semicomplete composition.

The following result can be found in the literature, see e.g. Theorem~2.2.9 of \cite{Bang-Jensen-Havet}.
%\begin{pro}\label{pro01}
%Every strong semicomplete digraph on $n \geq 3$ vertices contains a strong spanning tournament.
%\end{pro}
\begin{thm}\label{thm02} Every strong semicomplete digraph is vertex-pancyclic.
\end{thm}

It is well known that a tournament $T$ contains an $x-y$ Hamiltonian path if and only if there is an acyclic spanning subgraph $R$ (not necessarily induced) such that for each vertex $z$ of $T$, $R$ contains an $x-z$ path and a $z-y$ path \cite{Thomassen2}. Bang-Jensen and Huang \cite{Bang-Jensen-Huang1} proved that if a quasi-transitive digraph has both in- and out-branchings then it always contains such an acyclic spanning subdigraph. For semicomplete compositions, we give the following sufficient condition to guarantee the existence of this type of subdigraph.

\begin{thm}\label{thmf}
Let $Q=T[H_1, \dots, H_t]$ be a semicomplete composition. Then it contains an acyclic spanning subdigraph $R$ with a source $x$ and a sink $y$ such that for each vertex $z$ of $Q$, $R$ contains an $x-z$ path and a $z-y$ path, if one of the following assertions holds:
\begin{description}
\item[(a)] $Q$ is non-strong and has both in- and out-branchings.
\item[(b)] $Q$ is strong with $|V(H_i)|\geq 2$ for each $i\in [t]$.
\end{description}
\end{thm}
\begin{pf}
We first assume that $(a)$ holds. Let $Q$ have an out-branching, say $B^+_x$, rooted at $x$ and an in-branching, $B^-_y$, rooted at $y$, then $Q$ has precisely one initial strong component, say $Q'$, and precisely one strong terminal component, say $Q''$. Furthermore, we must have $x\in V(Q')$ and $y\in V(Q'')$. By the definition of a semicomplete composition, $Q'$ ($Q''$) is a subdigraph of some $H_i$ when $Q'$ ($Q''$) contains vertices from only one $H_i$, or is the union of some $H_i$ when $Q'$ ($Q''$) contains vertices from at least two $H_i$s.

We just consider the case that $Q'$ is a subdigraph of $H_1$, and $Q''=\bigcup_{i=s}^t{H_i}$ since the other cases are similar. Since $Q'$ (resp. $Q''$) is the unique strong initial (resp. terminal) strong component, we have $V(H_1)\Rightarrow \bigcup_{i=2}^{s-1}{H_i}\Rightarrow \bigcup_{i=s}^t{H_i}$. Let $B'^+_x$ be the subdigraph of $B^+_x$ induced by $V(H_1)$, it is not hard to see that $B'^+_x$ is an out-branching of $H_1$ rooted at $x$. Since $Q''$ is strong, it has an in-branching $B'^-_y$ rooted at $y$. We now construct a subdigraph $R$ of $Q$ from $B'^+_x$ and $B'^-_y$ by adding all arcs from $V(H_1)$ to $v$ and all arcs from $v$ to $V(Q'')$ for each $v\in \bigcup_{i=2}^{s-1}{V(H_i)}$ (if $\bigcup_{i=2}^{s-1}{V(H_i)}=\emptyset$, then we just add all arcs from $V(H_1)$ to $V(Q'')$). It can be checked that $R$ is the desired acyclic spanning subdigraph.

We next assume that $(b)$ holds. By Observation~\ref{obs1}, $T$ is strong. Furthermore, $T$ has a Hamiltonian cycle $C: u_1, u_2, \dots, u_t, u_1$ by Theorem~\ref{thm02}.

If $t\geq 3$, we construct a digraph $R$ from $C$ as follows: For each $i\in [t]$, substitute a copy of $H_i$ for $u_i$ and then delete all arcs inside $H_i$. Let $x\in V(H_1)$ and $y\in V(H_2)$; delete all arcs from $V(H_t)$ to $x$, all arcs from $V(H_1)\setminus \{x\}$ to $V(H_2)\setminus \{y\}$ and all arcs from $y$ to $V(H_3)$. It is not hard to check that $R$ is an acyclic spanning subdigraph of $Q$ satisfying the desired properties.

For the case that $t=2$, we construct a digraph $R$ from $C$ as follows: For each $i\in [2]$, substitute a copy of $H_i$ for $u_i$ and then delete all arcs inside $H_i$. Let $x\in V(H_1)$ and $y\in V(H_2)$; delete all arcs from $V(H_1)\setminus \{x\}$ to $V(H_2)\setminus \{y\}$, from $V(H_2)$ to $x$, and all arcs from $y$ to $V(H_1)$. Observe that $R$ is the desired acyclic spanning subdigraph of $Q$.
\end{pf}
%It is known that every strong digraph $D$ has an out- and in-branching rooted at any vertex of $D$.

Recall that we use the Hamiltonicity of a strong semicomplete digraph in the proof for the case that $Q$ is strong with $|V(H_i)|\geq 2$ for each $i\in [t]$ in Theorem~\ref{thmf}. In fact, the proof also means that the following more general result holds.

\begin{thm}\label{thmh}
Let $Q=T[H_1, \dots, H_t]$ be a digraph composition. If $T$ is Hamiltonian and $|V(H_i)|\geq 2$ for each $i\in [t]$, then $Q$ contains an acyclic spanning subdigraph $R$ with a source $x$ and a sink $y$ such that for each vertex $z$ of $Q$, $R$ contains an $x-z$ path and a $z-y$ path.
\end{thm}

\vskip 1cm

\noindent {\bf Acknowledgement.} Yuefang Sun was supported by Zhejiang Provincial Natural Science Foundation (No. LY20A010013) and National Natural Science Foundation of China (No. 11401389).

%\end{CJK}

\begin{thebibliography} {99}

\bibitem{Aho-Garey-Ullman}
A.V. Aho, M.R. Garey and J.D. Ullman, The transitive reduction of a directed graph, SIAM J. Comput. 1(2), 1972, 131--137.

%\bibitem{Aouchiche-Hansen}
%M. Aouchiche, P. Hansen, A survey of Nordhaus-Gaddum type relations,
%Discrete Appl. Math. 161(4/5), 2013, 466--546.

%\bibitem{Bang-Jensen} J. Bang-Jensen, Locally semicomplete digraphs: a generalization of tournaments. J. Graph Theory 14, 1990, 371--390.

\bibitem{Bang-Jensen-Gutin}
J. Bang-Jensen and G. Gutin, Digraphs: Theory, Algorithms and Applications, 2nd Edition, Springer, London, 2009.

%\bibitem{Bang-Jensen-Gutin2} J. Bang-Jensen and G. Gutin, Basic Terminology, Notation and Results, in {\em Classes of Directed Graphs}(J. Bang-Jensen and G. Gutin, eds.), Springer, 2018.

\bibitem{Bang-Jensen-Gutin-Huang}
J. Bang-Jensen, G. Gutin and J. Huang, A sufficient condition for a complete multipartite digraph to be hamiltonian, Discrete Math. 161, 1996, 1--12.

\bibitem{Bang-Jensen-Gutin-Yeo} J. Bang-Jensen, G. Gutin and A. Yeo,
Arc-disjoint strong spanning subgraphs of semicomplete compositions, J. Graph Theory DOI: 10.1002/jgt.22568.

\bibitem{Bang-Jensen-Havet} J. Bang-Jensen and F. Havet, Tournaments and Semicomplete Digraphs, in {\em Classes of Directed Graphs}
(J. Bang-Jensen and G. Gutin, eds.), Springer, 2018.

\bibitem{Bang-Jensen-Huang1} J. Bang-Jensen and J. Huang, Quasi-transitive digraphs, J. Graph Theory 20(2), 1995, 141--161.

%\bibitem{Bang-Jensen-Huang2} J. Bang-Jensen and J. Huang, Kings in quasi-transitive digraphs, Discrete Math. 185, 1998, 19--27.

%\bibitem{Bang-Jensen-Huang} J. Bang-Jensen and J. Huang, Decomposing locally semicomplete digraphs into strong spanning subdigraphs, J. Combin. Theory Ser. B, 102, 2012, 701--714.

\bibitem{Bang-Jensen-Huang-Yeo}
J. Bang-Jensen, J. Huang and A. Yeo, Strongly connected spanning subgraphs with the minimum number of arcs in quasi-transitive digraphs, SIAM J. Discrete Math. 16, 2003, 335--343.

%\bibitem{Bang-Jensen-Jordan} J. Bang-Jensen and T. Jord\'{a}n, Spanning 2-strong tournaments in 3-strong semicomplete digraphs, Discrete Math. 310, 2010, 1424--1428.

%\bibitem{Bang-Jensen-Kriesell} J. Bang-Jensen and M. Kriesell, Disjoint sub(di)graphs in digraphs, Electron. Notes Discrete Math. 34, 2009, 179--183.

\bibitem{Bang-Jensen-Nielsen-Yeo}
J. Bang-Jensen, M.H. Nielsen, A. Yeo, Longest path partitions in generalizations of tournaments, Discrete Math. 306(16), 2006, 1830--1839.

%\bibitem{Bang-Jensen-Yeo}J. Bang-Jensen and A. Yeo, Decomposing $k$-arc-strong tournaments into strong spanning subdigraphs, Combinatorica 24(3), 2004, 331--349.

%\bibitem{Bang-Jensen-Yeo2} J. Bang-Jensen and A. Yeo, Arc-disjoint spanning sub(di)graphs in Digraphs, Theoret. Comput. Sci. 438, 2012, 48--54.

%\bibitem{Bang-Jensen-Larsen-Maddaloni} J. Bang-Jensen, T.M. Larsen and A. Maddaloni, Disjoint paths in decomposable digraphs, J. Graph Th. 85(2) 2017, 545--567.

\bibitem{Bondy} J.A. Bondy and U.S.R. Murty, Graph Theory, Springer, Berlin, 2008.

%\bibitem{Camion} P. Camion, Chemins et circuits hamiltoniens des graphes complets, Comptes Rendus de l'Acad{\'e}mie des Sciences de Paris, 249, 1959, 2151-- 2152.

\bibitem{GSHC} H. Galeana-S\'anchez and C. Hern\'andez-Cruz, Quasi-transitive Digraphs and Their Extensions, in {\em Classes of Directed Graphs} (J. Bang-Jensen and G. Gutin, eds.), Springer, 2018.

%\bibitem{BT}  F. Boesch and R. Tindel, Robbins' theorem for mixed
%multigraphs.  Amer. Math. Monthly 87, 1980, 716--719.

%\bibitem{Chen-Li-Liu-Mao} L. Chen, X. Li, M. Liu and Y. Mao,  A solution to a conjecture on the generalized connectivity of graphs, J. Combin. Optim. 33(1), 2017, 275--282.

%\bibitem{Cheriyan-Salavatipour} J. Cheriyan and M. Salavatipour, Hardness and approximation results for packing Steiner trees, Algorithmica, 45, 2006, 21--43.

%\bibitem{Chud-Scott-SeymourAM} M. Chudnovsky, A. Scott and P.D. Seymour. Disjoint paths in tournaments. Adv. Math., 270, 2015, 582--597.

%\bibitem{Chud-Scott-Seymour} M. Chudnovsky, A. Scott and P.D. Seymour. Disjoint paths in unions of tournaments. arXiv:1604.02317, April 2016.

%\bibitem{DeVos-McDonald-Pivotto} M. DeVos, J. McDonald, I. Pivotto, Packing Steiner trees, J. Combin. Theory Ser. B, 119, 2016, 178--213.

%\bibitem{Edmonds} J. Edmonds, Edge-disjoint branchings, in {\em Combinatorial Algorithms} (B. Rustin ed.), Academic Press, 1973, 91--96.

%\bibitem{Feige-Halldorsson-Kortsarz-Srinivasan} U. Feige, M. Halldorsson, G. Kortsarz, and A. Srinivasan, Approximating the domatic number, SIAM J. Comput. 32(1), 2002, 172--195.

%\bibitem{Fortune-Hopcroft-Wyl} S. Fortune, J. Hopcroft and J. Wyllie, The directed subgraphs homeomorphism problem, Theoret. Comput. Sci. 10, 1980, 111--121.

%\bibitem{Guo} Y. Guo, Spanning local tournaments in locally semicomplete digraphs, Discrete Appl. Math. 79(1-3), 1997, 119--125.

\bibitem{Gutin2} G. Gutin, Cycles and paths in directed graphs, PhD thesis,
School of Mathematics, Tel Aviv University, 1993.

%\bibitem{Gutin3} G. Gutin, Finding a longest path in a complete multipartite digraph, SIAM J. Discrete Math. 6, 270--273, 1993.

\bibitem{Gutin4}
G. Gutin, Polynomial algorithms for finding Hamiltonian paths and cycles in
quasi-transitive digraphs, Australas. J. Combin. 10, 1994, 231--236.

\bibitem{Gutin} G. Gutin, Characterizations of vertex pancyclic and pancyclic
ordinary of complete multipartite digraphs, Discrete Math. 141(1-3), 1995, 153--162.

\bibitem{Gutin-Sun} G. Gutin, Y. Sun, Arc-disjoint in- and out-branchings rooted at the same vertex in compositions of digraphs, Discrete Math. 343(5), 2020, 111816.

%\bibitem{Gutin-Yeo} G. Gutin, A. Yeo, Kings in semicomplete multipartite digraphs, J. Graph Theory 33, 2000, 177--183.

%\bibitem{Dalmazzo} M. Dalmazzo, Nombre d'arcs dans les graphes $k$-arc-fortement connexes minimaux. C.R. Acad. Sci. Paris A 2853, 1977, 341--344.

%\bibitem{Hager} M. Hager, Pendant tree-connectivity, J. Combin. Theory Ser. B 38, 1985, 179--189.

\bibitem{Hammack} R.H. Hammack, Digraph Products, in {\em Classes of Directed Graphs} (J. Bang-Jensen and G. Gutin, eds.), Springer, 2018.

%\bibitem{Huang-Li} J. Huang, W. Li, Toppling kings in a tournament by introducing new kings, J. Graph Theory 11, 1987, 7--11.

%\bibitem{JRST} T. Johnson, N. Robertson, P.D. Seymour and R. Thomas, Directed Tree-Width, J. Combin. Th. Ser. B 82(1), 2001, 138--154.

%\bibitem{Kriesell} M. Kriesell, Edge-disjoint trees containing some given vertices in a graph, J. Combin. Theory Ser. B 88, 2003, 53--65.

%\bibitem{Koh-Tan} K.M. Koh and B.P. Tan, Number of 4-kings in bipartite tournaments with no 3-kings, Discrete Math. 154(1--3), 1996, 281--287.

%\bibitem{Koh-Tan2} K.M. Koh and B.P. Tan, The number of kings in a multipartite tournament, Discrete Math. 167/168, 1997, 411--418.

\bibitem{Khuller-Raghavachari-Young}
S. Khuller, B. Raghavachari and N. Young, Approximating the minimum equivalent digraph, SIAM J. Comput. 24, 1995, 859--872.

%\bibitem{Landau} H.G. Landau, On dominance relations and the structure of animal societies III, The condition for a score structure, Bull. Math. Biophys. 15, 1953, 143--148.

%\bibitem{Lau} L. Lau, An approximate max-Steiner-tree-packing min-Steiner-cut theorem, Combinatorica 27, 2007, 71--90.

%\bibitem{LiThesis} S. Li, Some topics on generalized connectivity of graphs, PhD thesis, Nankai University, 2012.

%\bibitem{Li-Li}
%S. Li and X. Li, Note on the hardness of generalized connectivity,
%J. Comb. Optim. 24(3), 2012, 389--396.

%\bibitem{Li-Li-Zhou} S. Li, X. Li and W. Zhou, Sharp bounds for the generalized connectivity $\kappa_3(G)$, Discrete Math. 310, 2010, 2147--2163.

%\bibitem{Li-Mao} X. Li and Y. Mao, A survey on the generalized connectivity of graphs, arXiv:1207.1838, v10, Aug 2015.

%\bibitem{Li-Mao5} X. Li and Y. Mao, Generalized Connectivity of Graphs, Springer, Switzerland, 2016.

%\bibitem{Li-Mao-Sun}
%X. Li, Y. Mao and Y. Sun, On the generalized (edge-)connectivity of
%graphs, Australas. J. Combin. 58(2), 2014, 304--319.

%\bibitem{Mader} W. Mader, Minimal $n$-fach zusammenh\"{a}ngende Digraphen, J. Combin.
%Theory Ser. B 38(2), 1985, 102--117.

%\bibitem{Nash-Williams}
%C.St.J.A. Nash-Williams, Edge-disjonint spanning trees of finite
%graphs, J. London Math. Soc. 36, 1961, 445--450.

%\bibitem{Ng} L.L. Ng. Hamiltonian decomposition of lexicographic products of digraphs. J. Combin. Theory Ser. B, 73(2), 1998, 119--129.

%\bibitem{Ozeki-Yamashita}
%K. Ozeki, T. Yamashita, Spanning trees: a survey, Graphs Combin.
%27(1), 2011, 1--26.

%\bibitem{Palmer}
%E. Palmer, On the spanning tree packing number of a graph: a survey,
%Discrete Math. 230, 2001, 13--21.

%\bibitem{Reid} K.B. Reid, Every vertex a king, Discrete Math. 38, 1982, 93--98.

%\bibitem{Schr} A. Schrijver, Finding $k$ partially disjoint paths in a directed planar graph. SIAM J. Comput. 23(4), 1994, 780--788.

\bibitem{Sun} Y. Sun, Kings in compositions of digraphs, in preparation.

%\bibitem{Sun-Gutin} Y. Sun, G. Gutin, Strong subgraph $k$-connectivity bounds, arXiv:1803.00281v1 [cs.DM] 1 Mar 2018.

%\bibitem{Sun-Gutin2} Y. Sun, G. Gutin, Strong subgraph connectivity of digraphs, submitted.

%\bibitem{Sun-Gutin3} Y. Sun, G. Gutin, Strong subgraph connectivity of digraphs: a survey, arXiv:1808.02740v1 [cs.DM] 8 Aug 2018.

\bibitem{Sun-Gutin-Ai} Y. Sun, G. Gutin, J. Ai, Arc-disjoint strong spanning subdigraphs in compositions and products of digraphs, Discrete Math. 342(8), 2019, 2297--2305.

%\bibitem{Sun-Gutin-Yeo-Zhang} Y. Sun, G. Gutin, A. Yeo, X. Zhang, Strong subgraph $k$-connectivity, J. Graph Theory 92(1), 2019, 5--18.

%\bibitem{Sun-Zhang} Y. Sun, X. Zhang, Complexity and computation of strong subgraph packing in digraphs, in preparation.





%\bibitem{shiloachIPL8} Y. Shiloach, Edge-disjoint branching in directed multigraphs, Inf. Process. Lett. 8(1), 1979, 24--27.


%\bibitem{Sun-Li}
%Y. Sun and X. Li, On the difference of two generalized connectivities
%of a graph, J. Comb. Optim. 33(1), 2017, 283--291.

%\bibitem{Thom} C. Thomassen, Highly connected non-2-linked digraphs, Combinatorica 11(4), 1991, 393--395.

%\bibitem{Thomassen}C. Thomassen, Configurations in graphs of large minimum degree, connectivity, or chromatic number, Annals of the New York Academy of Sciences, 555, 1989, 402--412.

\bibitem{Thomassen2}
C. Thomassen, Hamiltonian-connected tournaments, J. Combin. Theory Ser. B 28, 1980, 142--163.

%\bibitem{Thomassen1} C. Thomassen, Connectivity in tournaments, Graph Theory and Combinatorics, (B. Bollodis, ed.), Academic Press, London, 1984, 305--313.

%\bibitem{Tillson}
%T.W. Tillson, A Hamiltonian decomposition of $K^*_{2m}$, $2m \geq
%8$, J. Combin. Theory Ser. B 29(1), 1980, 68--74.

%\bibitem{Trotter-Erdos}
%W.T. Trotter, Jr. and P. Erd\H{o}s. When the Cartesian product of
%directed cycles is Hamiltonian, J. Graph Theory, 2(2), 1978, 137--142.

%\bibitem{Tutte}
%W. Tutte, On the problem of decomposing a graph into $n$ connected
%factors, J. London Math. Soc. 36, 1961, 221--230.

%\bibitem{Volkmann} L. Volkmann, Cycles in multipartite tournaments: results and problems, Discrete Math. 245(1), 2002, 19--53.

%\bibitem{Wang-Zhang} R. Wang and H. Zhang, $(k+1)$-kernels and the number of $k$-kings in $k$-quasi-transitive digraphs, Discrete Math. 338(1), 2015, 114--121.

%\bibitem{West-Wu} D. West, H. Wu, Packing Steiner trees and S-connectors in graphs, J. Combin. Theory Ser. B 102, 2012, 186--205.
\end{thebibliography}
\end{document}